\documentclass [12pt]{article}
\usepackage{graphicx,amssymb,amsfonts,latexsym,amsmath,amsthm,times}
\usepackage{epsfig}
\usepackage{color}
\usepackage[english]{babel}
\usepackage[a4paper]{geometry}
\def\lanbox{\hbox{$\, \vrule height 0.25cm width 0.25cm depth 0.01cm \,$}}

\numberwithin{equation}{section}

\begin{document}

\vspace*{1.4cm}

\normalsize \centerline{\Large \bf ON THE WELL-POSEDNESS
OF SOME MODEL}

\medskip

\centerline{\Large \bf ARISING IN THE MATHEMATICAL BIOLOGY}

\vspace*{1cm}

\centerline{\bf Messoud Efendiev$^{1,2}$, Vitali Vougalter$^{3 \ *}$}

\bigskip

\centerline{$^1$ Helmholtz Zentrum M\"unchen, Institut f\"ur Computational
Biology, Ingolst\"adter Landstrasse 1}

\centerline{Neuherberg, 85764, Germany}

\centerline{e-mail: messoud.efendiyev@helmholtz-muenchen.de}

\centerline{$^2$ Department of Mathematics, Marmara University, Istanbul,
T\"urkiye}

\centerline{e-mail: m.efendiyev@marmara.edu.tr}

\bigskip

\centerline{$^{3 \ *}$  Department of Mathematics, University
of Toronto}

\centerline{Toronto, Ontario, M5S 2E4, Canada}

\centerline{ e-mail: vitali@math.toronto.edu}

\medskip


\vspace*{0.25cm}

\noindent {\bf Abstract:}
In the article we establish the global well-posedness in
$W^{1,2,2}({\mathbb R}\times {\mathbb R}^{+})$ of the
integro-differential equation in the case of anomalous diffusion when the
one dimensional negative Laplace operator is raised to a fractional power in
the presence of the transport term. The model is relevant to the cell
population dynamics in the Mathematical Biology.
Our proof relies on a fixed point technique.

\vspace*{0.25cm}



\vspace*{0.5cm}

\bigskip

\bigskip


\setcounter{section}{1}

\centerline{\bf 1. Introduction}

\medskip

\noindent
The present work deals with the global well-posedness
of the nonlocal reaction-diffusion equation with the constants
$\displaystyle{0<\alpha<1, \ a\geq 0}$ and $b\in {\mathbb R}$, namely
\begin{equation}
\label{h}
\frac{\partial u}{\partial t} =
-\Bigg(-\frac{\partial^{2}}{\partial x^{2}}\Bigg)^{\alpha} u +
b\frac{\partial u}{\partial x}+au+
\int_{-\infty}^{\infty}G(x-y)F(u(y,t), y)dy, \quad
\end{equation}
which appears in the cell population dynamics. Let us assume that the initial
condition for (\ref{h}) is given by
\begin{equation}
\label{ic}
u(x,0)=u_{0}(x)\in H^{2}({\mathbb R}).
\end{equation}  
Note that the existence of stationary solutions of the equation analogous to
(\ref{h}) with $\displaystyle{\alpha=\frac{1}{2}}$ on the whole real line and
on a finite interval with periodic boundary conditions was discussed in
~\cite{EV22}. The case of the normal diffusion and transport was treated in
~\cite{EV20}. The article ~\cite{VV21} is devoted to the existence of
stationary solutions of the integro-differential equation involving the
fractional Laplacian in one dimension raised to the power
$\displaystyle{0<\alpha<\frac{1}{4}}$, the drift and the cell influx/efflux.
The situation without the transport term was considered
in ~\cite{VV18}. Spatial structures and generalized travelling waves for an
integro-differential equation were covered in ~\cite{ABVV10}. In work
~\cite{VV130} the authors deal with the 
emergence and propagation of patterns in nonlocal reaction-
diffusion equations arising in the theory of speciation and involving the
drift term. The existence of steady states and travelling waves for the
non-local Fisher-KPP equation was established in ~\cite{BNPR09}. In
~\cite{BHN05} the authors estimated the speed of propagation for KPP type
problems in the periodic framework. Important applications to the theory of
reaction-diffusion equations with non-Fredholm operators were developed in
~\cite{DMV05}, ~\cite{DMV08}.  

\noindent
The space variable $x$ in our article corresponds to the cell genotype,
$u(x,t)$ stands for the cell density as a function of the genotype and time.
The right side of (\ref{h}) describes the evolution of the cell density due to
the cell proliferation, mutations and transport.
The anomalous diffusion term here is correspondent to the change of genotype
via the small random mutations, and the integral term describes large mutations.
The function $F(u, x)$ denotes the rate of cell birth dependent on $u$ and $x$
(density dependent proliferation), and the kernel $G(x-y)$ gives
the proportion of newly born cells changing their genotype from $y$ to $x$.
We assume that it depends on the distance between the genotypes.

\noindent
The operator $\displaystyle{\Bigg(-\frac{\partial^{2}}{\partial x^{2}}
\Bigg)^{\alpha}}$ in problem (\ref{h}) describes a particular case of the
anomalous diffusion actively studied in the context of different applications
in: plasma physics and
turbulence \cite{2}, \cite{1}, surface diffusion \cite{3}, \cite{4},
semiconductors \cite{5} and so on. It is actively used in the works on the
nonlocal diffusive processes. The probabilistic realization of the
anomalous diffusion was covered in ~\cite{MK00}. Front propagation equations
with anomalous diffusion were studied extensively in
recent years (see e.g. ~\cite{VNN10}, ~\cite{VNN13}).
The operator
$\displaystyle{\Bigg(-\frac{\partial^{2}}{\partial x^{2}}\Bigg)^{\alpha}}$
is defined by virtue of the spectral calculus. This is the pseudo-differential
operator with symbol $\displaystyle{|p|^{2\alpha}}$, such that
$$
\Bigg(-\frac{d^{2}}{dx^{2}}\Bigg)^{\alpha}\phi(x)=\frac{1}{\sqrt{2\pi}}
\int_{-\infty}^{\infty}|p|^{2\alpha}{\widehat \phi}(p)e^{ipx}dp, \quad
\phi(x)\in H^{2\alpha}({\mathbb R}),
$$
where the Sobolev space is defined in (\ref{ss}).
The standard Fourier transform in this context is given by
\begin{equation}
\label{ft}
\widehat{\phi}(p)=\frac{1}{\sqrt{2\pi}}\int_{-\infty}^{\infty}\phi(x)e^{-ipx}dx,
\quad p\in {\mathbb R}.
\end{equation}
Evidently, the upper bound
\begin{equation}
\label{fub}  
\|\widehat{\phi}(p)\|_{L^{\infty}({\mathbb R})}\leq \frac{1}{\sqrt{2\pi}}
\|\phi(x)\|_{L^{1}({\mathbb R})}
\end{equation}
holds (see e.g. ~\cite{LL97}). Clearly, (\ref{fub}) yields
\begin{equation}
\label{fub1}  
\|p^{2}\widehat{\phi}(p)\|_{L^{\infty}({\mathbb R})}\leq \frac{1}{\sqrt{2\pi}}
\Big\|\frac{d^{2}\phi}{dx^{2}}\Big\|_{L^{1}({\mathbb R})}.
\end{equation}
We suppose that the conditions below on the integral kernel involved in 
problem (\ref{h}) are fulfilled.

\medskip

\noindent
{\bf Assumption 1.1.}  {\it Let $G(x): {\mathbb R}\to {\mathbb R}$ be
nontrivial, so that
$\displaystyle{G(x), \ \frac{d^{2}G(x)}{dx^{2}}\in L^{1}({\mathbb R})}$.}     

\medskip

\noindent
This enables us to introduce the technical quantity
\begin{equation}
\label{q}  
Q:=\sqrt{\|G(x)\|_{L^{1}({\mathbb R})}^{2}+
\Big\|\frac{d^{2}G}{dx^{2}}\Big\|_{L^{1}({\mathbb R})}^{2}}.
\end{equation}
Hence, $0<Q<\infty$.

\medskip

\noindent
From the perspective of the applications, the space dimension is
not restricted to $d=1$ because the space variable is correspondent to the cell
genotype but not to the usual physical space.
We have the Sobolev space
\begin{equation}
\label{ss}  
H^{2\alpha}({\mathbb R}):=
\Bigg\{\phi(x):{\mathbb R}\to {\mathbb {\mathbb R}} \ | \
\phi(x)\in L^{2}({\mathbb R}), \ \Bigg(-\frac{d^{2}}{dx^{2}}\Bigg)^{\alpha}\phi \in
L^{2}({\mathbb R}) \Bigg \}, 
\end{equation}
where $0<\alpha\leq 1$.
It is equipped with the norm
\begin{equation}
\label{n}
\|\phi\|_{H^{2\alpha}({\mathbb R})}^{2}:=\|\phi\|_{L^{2}({\mathbb R})}^{2}+
\Bigg\| \Bigg(-\frac{d^{2}}{dx^{2}}\Bigg)^{\alpha}\phi\Bigg\|_{L^{2}({\mathbb R})}^{2}.
\end{equation}
Obviously, in the particular case of $\displaystyle{\alpha=1}$, we have
\begin{equation}
\label{n1}
\|\phi\|_{H^{2}({\mathbb R})}^{2}:=\|\phi\|_{L^{2}({\mathbb R})}^{2}+
\Bigg\| \frac{d^{2}\phi}{dx^{2}}\Bigg\|_{L^{2}({\mathbb R})}^{2}.
\end{equation}
To establish the global well-posedness for problem (\ref{h}), (\ref{ic})
we will use the function space
$$
W^{1, 2, 2}({\mathbb R}\times [0, T]):=
$$
\begin{equation}
\label{122}
\Big\{u(x,t): {\mathbb R}\times [0, T]\to {\mathbb R} \ \Big|
\ u(x,t), \ \frac{\partial^{2}u}{\partial x^{2}}, \
\frac{\partial u}{\partial t}
\in L^{2}({\mathbb R}\times [0, T]) \Big\},
\end{equation}  
such that
$$
\|u(x,t)\|_{W^{1, 2, 2}({\mathbb R}\times [0, T])}^{2}:=
$$
\begin{equation}
\label{122n}
\Big\|\frac{\partial u}{\partial t}\Big\|_{L^{2}({\mathbb R}\times [0, T])}^{2}+
\Big\|\frac{\partial^{2}u}{\partial x^{2}}\Big\|_{L^{2}({\mathbb R}\times [0, T])}^{2}+
\|u\|_{L^{2}({\mathbb R}\times [0, T])}^{2}
\end{equation}  
with $T>0$. Here
$$
\|u\|_{L^{2}({\mathbb R}\times [0, T])}^{2}:=\int_{0}^{T}\int_{-\infty}^{\infty}
|u(x,t)|^{2}dxdt.
$$
Throughout the article we will also use another norm
$$
\|u(x,t)\|_{L^{2}({\mathbb R})}^{2}:=\int_{-\infty}^{\infty}|u(x,t)|^{2}dx.
$$

\medskip

\noindent
{\bf Assumption 1.2.} {\it Function
$F(u, x): {\mathbb R}\times{\mathbb R}\to {\mathbb R}$ is satisfying the
Caratheodory condition (see ~\cite{K64}), such that
\begin{equation}
\label{lub}
|F(u, x)|\leq k|u|+h(x) \quad for \quad u\in {\mathbb R}, \quad x\in {\mathbb R}
\end{equation}  
with a constant $k>0$ and
$h(x): {\mathbb R}\to {\mathbb R}^{+}, \ h(x)\in L^{2}({\mathbb R})$.
Furthermore, it is a Lipschitz continuous function, so that
\begin{equation}
\label{lip}
|F(u_{1}, x)-F(u_{2}, x)|\leq l|u_{1}-u_{2}| \quad for \quad any \quad
u_{1, 2}\in {\mathbb R}, \quad x\in {\mathbb R}
\end{equation}
with a constant $l>0$.}

\medskip

\noindent
In the work ${\mathbb R}^{+}$ denotes the nonnegative semi-axis.
The solvability of a local elliptic problem in a bounded domain in
${\mathbb R}^{N}$ was discussed in ~\cite{BO86}. The nonlinear function there
was allowed to have a sublinear growth.
Let us apply the standard Fourier transform (\ref{ft}) to both sides of
problem (\ref{h}), (\ref{ic}). This yields
\begin{equation}
\label{hf}
\frac{\partial \widehat{u}}{\partial t}=[-|p|^{2\alpha}+ibp+a]\widehat{u}+
\sqrt{2\pi}\widehat{G}(p)\widehat{f}_{u}(p,t),
\end{equation}
\begin{equation}
\label{icf}
\widehat{u(x,0)}(p)=\widehat{u_{0}}(p).  
\end{equation}
Here and further down $\widehat{f}_{u}(p,t)$ will stand for the Fourier image
of $F(u(x,t), x)$. Clearly, we have
$$
u(x,t)=\frac{1}{\sqrt{2\pi}}\int_{-\infty}^{\infty}\widehat{u}(p,t)e^{ipx}dp,
\quad 
\frac{\partial u}{\partial t}=\frac{1}{\sqrt{2\pi}}\int_{-\infty}^{\infty}
\frac{\partial \widehat{u}(p,t)}{\partial t}e^{ipx}dp
$$
with $x\in {\mathbb R}, \ t\geq 0$. By means of the Duhamel's principle,
we can reformulate problem (\ref{hf}), (\ref{icf}) as
$$
\widehat{u}(p,t)=
$$
\begin{equation}
\label{duh}
e^{t\{-|p|^{2\alpha}+ibp+a\}}\widehat{u_{0}}(p)+\int_{0}^{t}
e^{(t-s)\{-|p|^{2\alpha}+ibp+a\}}\sqrt{2\pi}\widehat{G}(p)\widehat{f}_{u}(p,s)ds.
\end{equation}
Related to equation (\ref{duh}), we consider the auxiliary problem
$$
\widehat{u}(p,t)=
$$
\begin{equation}
\label{aux}
e^{t\{-|p|^{2\alpha}+ibp+a\}}\widehat{u_{0}}(p)+\int_{0}^{t}
e^{(t-s)\{-|p|^{2\alpha}+ibp+a\}}\sqrt{2\pi}\widehat{G}(p)\widehat{f}_{v}(p,s)ds,
\end{equation}
where $\widehat{f}_{v}(p,s)$ denotes the Fourier image of
$F(v(x,s), x)$ under transform (\ref{ft}).

\noindent
Let us introduce the operator $t_{a, b}$, such that $u =t_{a, b}v$, where $u$ 
satisfies (\ref{aux}). Our main result is as follows.

\bigskip

\noindent
{\bf Theorem 1.3.} {\it Let Assumptions 1.1 and 1.2 hold and
\begin{equation}
\label{qlt}
Ql\sqrt{T^{2}e^{2aT}(1+2[a+|b|+1]^{2})+1}<1.  
\end{equation}
Then equation (\ref{aux}) defines the map
$t_{a, b}: W^{1, 2, 2}({\mathbb R}\times [0, T])\to W^{1, 2, 2}
({\mathbb R}\times [0, T])$, which is a strict contraction.
The unique fixed point $w(x,t)$ of this map $t_{a, b}$ is the only solution of
problem (\ref{h}), (\ref{ic}) in
$W^{1, 2, 2}({\mathbb R}\times [0, T])$.}

\medskip

\noindent
The final statement of the article deals with the global well-posedness for
our equation. 

\medskip

\noindent
{\bf Corollary 1.4.} {\it Let the assumptions of Theorem 1.3 be valid.
Then problem (\ref{h}), (\ref{ic}) admits a unique solution   
$w(x,t)\in W^{1, 2, 2}({\mathbb R}\times {\mathbb R}^{+})$. This solution is
nontrivial for $x\in {\mathbb R}$ and $t\in {\mathbb R}^{+}$ 
provided the intersection of supports of the Fourier transforms of functions
$\hbox{supp}\widehat{F(0, x)}\cap \hbox{supp}\widehat{G}$ is a set of
nonzero Lebesgue measure in ${\mathbb R}$.}

\medskip

\noindent
We proceed to the proof of our main proposition.

\bigskip


\setcounter{section}{2}
\setcounter{equation}{0}

\centerline{\bf 2. The well-posedness of the model}

\bigskip

\noindent
{\it Proof of Theorem 1.3.} Let us choose an arbitrary
$v(x,t)\in W^{1, 2, 2}({\mathbb R}\times [0, T])$.

\noindent
It can be easily verified that the first term in the right side of (\ref{aux})
belongs to $L^{2}({\mathbb R}\times [0, T])$. Indeed,
$$
\|e^{t\{-|p|^{2\alpha}+ibp+a\}}\widehat{u_{0}}(p)\|_{L^{2}({\mathbb R})}^{2}=
\int_{-\infty}^{\infty}e^{-2t|p|^{2\alpha}}e^{2at}|\widehat{u_{0}}(p)|^{2}dp\leq
e^{2at}\|u_{0}\|_{L^{2}({\mathbb R})}^{2},
$$
such that
$$
\|e^{t\{-|p|^{2\alpha}+ibp+a\}}\widehat{u_{0}}(p)\|_{L^{2}({\mathbb R}\times [0, T])}^{2}=
\int_{0}^{T}\|e^{t\{-|p|^{2\alpha}+ibp+a\}}\widehat{u_{0}}(p)\|_{L^{2}({\mathbb R})}^{2}dt
\leq
$$
$$
\int_{0}^{T}e^{2at}\|u_{0}\|_{L^{2}({\mathbb R})}^{2}dt.
$$
Clearly, this equals to
$\displaystyle{\frac{e^{2aT}-1}{2a}\|u_{0}\|_{L^{2}({\mathbb R})}^{2}}$
for $a>0$ and
$T\|u_{0}\|_{L^{2}({\mathbb R})}^{2}$ if $a=0$. Hence,
\begin{equation}
\label{u0l2}  
e^{t\{-|p|^{2\alpha}+ibp+a\}}\widehat{u_{0}}(p)\in L^{2}({\mathbb R}\times [0, T]).
\end{equation}
We estimate the norm of the second term in the right side of (\ref{aux}) as
$$
\Big\|\int_{0}^{t}
e^{(t-s)\{-|p|^{2\alpha}+ibp+a\}}\sqrt{2\pi}\widehat{G}(p)\widehat{f}_{v}(p,s)ds\Big\|_
{L^{2}({\mathbb R})}\leq
$$
$$  
\int_{0}^{t}\Big\|
e^{(t-s)\{-|p|^{2\alpha}+ibp+a\}}\sqrt{2\pi}\widehat{G}(p)\widehat{f}_{v}(p,s)\Big\|_
{L^{2}({\mathbb R})}ds.
$$
Obviously,
$$  
\Big\|
e^{(t-s)\{-|p|^{2\alpha}+ibp+a\}}\sqrt{2\pi}\widehat{G}(p)\widehat{f}_{v}(p,s)ds\Big\|_
{L^{2}({\mathbb R})}^{2}=
$$
\begin{equation}
\label{intgfv}  
\int_{-\infty}^{\infty}e^{-2(t-s)|p|^{2\alpha}}e^{2a(t-s)}2\pi|\widehat{G}(p)|^{2}
|\widehat{f}_{v}(p,s)|^{2}dp.
\end{equation}
We use inequality (\ref{fub}) to obtain the upper bound on the right side of
(\ref{intgfv}) as
$$
e^{2a(t-s)}2\pi\|\widehat{G}(p)\|_{L^{\infty}({\mathbb R})}^{2}
\|F(v(x,s),x)\|_{L^{2}({\mathbb R})}^{2}\leq
$$
$$
e^{2aT}\|G(x)\|_{L^{1}({\mathbb R})}^{2}
\|F(v(x,s),x)\|_{L^{2}({\mathbb R})}^{2}.
$$
Thus,
$$  
\Big\|
e^{(t-s)\{-|p|^{2\alpha}+ibp+a\}}\sqrt{2\pi}\widehat{G}(p)\widehat{f}_{v}(p,s)ds\Big\|_
{L^{2}({\mathbb R})}\leq
$$
$$
e^{aT}\|G(x)\|_{L^{1}({\mathbb R})}
\|F(v(x,s),x)\|_{L^{2}({\mathbb R})}.
$$
By means of (\ref{lub}), we have
\begin{equation}
\label{lubs}  
\|F(v(x,s),x)\|_{L^{2}({\mathbb R})}\leq k\|v(x,s)\|_{L^{2}({\mathbb R})}+\|h(x)\|_
{L^{2}({\mathbb R})}.
\end{equation}
Therefore,
$$  
\Big\|
e^{(t-s)\{-|p|^{2\alpha}+ibp+a\}}\sqrt{2\pi}\widehat{G}(p)\widehat{f}_{v}(p,s)\Big\|_
{L^{2}({\mathbb R})}\leq
$$
$$
e^{aT}\|G(x)\|_{L^{1}({\mathbb R})}
\{k\|v(x,s)\|_{L^{2}({\mathbb R})}+\|h(x)\|_{L^{2}({\mathbb R})}\},
$$
such that
$$  
\Big\|\int_{0}^{t}
e^{(t-s)\{-|p|^{2\alpha}+ibp+a\}}\sqrt{2\pi}\widehat{G}(p)\widehat{f}_{v}(p,s)ds\Big\|_
{L^{2}({\mathbb R})}\leq
$$
$$
ke^{aT}\|G(x)\|_{L^{1}({\mathbb R})}\int_{0}^{T}\|v(x,s)\|_{L^{2}({\mathbb R})}ds+
Te^{aT}\|G(x)\|_{L^{1}({\mathbb R})}\|h(x)\|_{L^{2}({\mathbb R})}.
$$
By virtue of the Schwarz inequality
\begin{equation}
\label{sch}  
\int_{0}^{T}\|v(x,s)\|_{L^{2}({\mathbb R})}ds\leq
\sqrt{\int_{0}^{T}\|v(x,s)\|_{L^{2}({\mathbb R})}^{2}ds}\sqrt{T}.
\end{equation}
This implies that
$$
\Big\|\int_{0}^{t}
e^{(t-s)\{-|p|^{2\alpha}+ibp+a\}}\sqrt{2\pi}\widehat{G}(p)\widehat{f}_{v}(p,s)ds\Big\|_
{L^{2}({\mathbb R})}^{2}\leq
$$
$$
e^{2aT}\|G(x)\|_{L^{1}({\mathbb R})}^{2}
\{k\sqrt{T}\|v(x,s)\|_{L^{2}({\mathbb R}\times [0,T])}+T\|h(x)\|_{L^{2}({\mathbb R})}\}^{2}.
$$
Let us estimate the norm as
$$
\Big\|\int_{0}^{t}
e^{(t-s)\{-|p|^{2\alpha}+ibp+a\}}\sqrt{2\pi}\widehat{G}(p)\widehat{f}_{v}(p,s)ds\Big\|_
{L^{2}({\mathbb R}\times [0,T])}^{2}=
$$
$$
\int_{0}^{T}\Big\|\int_{0}^{t}
e^{(t-s)\{-|p|^{2\alpha}+ibp+a\}}\sqrt{2\pi}\widehat{G}(p)\widehat{f}_{v}(p,s)ds\Big\|_
{L^{2}({\mathbb R})}^{2}dt\leq 
$$
$$
e^{2aT}\|G(x)\|_{L^{1}({\mathbb R})}^{2}
\{k\|v(x,s)\|_{L^{2}({\mathbb R}\times [0,T])}+\sqrt{T}\|h(x)\|_{L^{2}({\mathbb R})}\}^{2}
T^{2}<\infty
$$
under the given conditions for
$v(x,s)\in W^{1, 2, 2}({\mathbb R}\times [0, T])$. Hence,
\begin{equation}
\label{0tetsgvl2}  
\int_{0}^{t}
e^{(t-s)\{-|p|^{2\alpha}+ibp+a\}}\sqrt{2\pi}\widehat{G}(p)\widehat{f}_{v}(p,s)ds\in
L^{2}({\mathbb R}\times [0,T]).
\end{equation}
By means of (\ref{u0l2}), (\ref{0tetsgvl2}) and (\ref{aux}), we arrive at
\begin{equation}
\label{upt12}  
\widehat{u}(p,t)\in L^{2}({\mathbb R}\times [0,T]).
\end{equation}
Therefore,
\begin{equation}
\label{uxtl2}  
u(x,t)\in L^{2}({\mathbb R}\times [0,T]).
\end{equation}
It follows easily from (\ref{aux}) that
$$
p^{2}\widehat{u}(p,t)=
$$
\begin{equation}
\label{aux2}
e^{t\{-|p|^{2\alpha}+ibp+a\}}p^{2}\widehat{u_{0}}(p)+\int_{0}^{t}
e^{(t-s)\{-|p|^{2\alpha}+ibp+a\}}\sqrt{2\pi}p^{2}\widehat{G}(p)\widehat{f}_{v}(p,s)ds.
\end{equation}
We consider the first term in the right side of (\ref{aux2}), such that
$$
\|e^{t\{-|p|^{2\alpha}+ibp+a\}}p^{2}\widehat{u_{0}}(p)\|_{L^{2}({\mathbb R}\times [0,T])}^{2}=
\int_{0}^{T}\int_{-\infty}^{\infty}e^{-2t|p|^{2\alpha}}e^{2at}|p^{2}\widehat{u_{0}}(p)|^{2}dp
dt\leq
$$
$$
\int_{0}^{T}\int_{-\infty}^{\infty}e^{2at}|p^{2}\widehat{u_{0}}(p)|^{2}dpdt.
$$
Evidently, this is equal to
$\displaystyle
{\frac{e^{2aT}-1}{2a}\Big\|\frac{d^{2}u_{0}}{dx^{2}}\Big\|_{L^{2}({\mathbb R})}^{2}}$ if
$a>0$ and
$\displaystyle{T\Big\|\frac{d^{2}u_{0}}{dx^{2}}\Big\|_{L^{2}({\mathbb R})}^{2}}$
for $a=0$. This means that
\begin{equation}
\label{p2u0hpl2}
e^{t\{-|p|^{2\alpha}+ibp+a\}}p^{2}\widehat{u_{0}}(p)\in L^{2}({\mathbb R}\times [0,T]).
\end{equation}  
Let us turn our attention to analyzing the second term in the right side of
(\ref{aux2}). Clearly,
$$
\Big\|\int_{0}^{t}
e^{(t-s)\{-|p|^{2\alpha}+ibp+a\}}\sqrt{2\pi}p^{2}\widehat{G}(p)\widehat{f}_{v}(p,s)ds
\Big\|_{L^{2}({\mathbb R})}\leq 
$$
$$
\int_{0}^{t}\Big\|
e^{(t-s)\{-|p|^{2\alpha}+ibp+a\}}\sqrt{2\pi}p^{2}\widehat{G}(p)\widehat{f}_{v}(p,s)
\Big\|_{L^{2}({\mathbb R})}ds.
$$
We have
$$
\Big\|e^{(t-s)\{-|p|^{2\alpha}+ibp+a\}}\sqrt{2\pi}p^{2}\widehat{G}(p)\widehat{f}_{v}(p,s)
\Big\|_{L^{2}({\mathbb R})}^{2}=
$$
\begin{equation}
\label{etsp2gpl2}
\int_{-\infty}^{\infty}
e^{-2(t-s)|p|^{2\alpha}}e^{2a(t-s)}2\pi|p^{2}\widehat{G}(p)|^{2}|\widehat{f}_{v}(p,s)|^{2}
dp.
\end{equation}  
The right side of (\ref{etsp2gpl2}) can be bounded above using (\ref{fub1})
as
$$
2\pi e^{2aT}\|p^{2}\widehat{G}(p)\|_{L^{\infty}({\mathbb R})}^{2}
\int_{-\infty}^{\infty}|\widehat{f}_{v}(p,s)|^{2}dp\leq
$$
$$
e^{2aT}\Big\|\frac{d^{2}G}{dx^{2}}\Big\|_{L^{1}({\mathbb R})}^{2}
\|F(v(x,s),x)\|_{L^{2}({\mathbb R})}^{2}.
$$
Let us use inequality (\ref{lubs}) to obtain the estimate
$$
\Big\|e^{(t-s)\{-|p|^{2\alpha}+ibp+a\}}\sqrt{2\pi}p^{2}\widehat{G}(p)\widehat{f}_{v}(p,s)
\Big\|_{L^{2}({\mathbb R})}\leq
$$
$$
e^{aT}\Big\|\frac{d^{2}G}{dx^{2}}\Big\|_{L^{1}({\mathbb R})}\{k\|v(x,s)\|_
{L^{2}({\mathbb R})}+\|h(x)\|_{L^{2}({\mathbb R})}\},
$$
such that
$$
\Big\|\int_{0}^{t}e^{(t-s)\{-|p|^{2\alpha}+ibp+a\}}\sqrt{2\pi}p^{2}\widehat{G}(p)
\widehat{f}_{v}(p,s)ds\Big\|_{L^{2}({\mathbb R})}\leq
$$
$$
ke^{aT}\Big\|\frac{d^{2}G}{dx^{2}}\Big\|_{L^{1}({\mathbb R})}\int_{0}^{T}\|v(x,s)\|_
{L^{2}({\mathbb R})}ds+
Te^{aT}\Big\|\frac{d^{2}G}{dx^{2}}\Big\|_{L^{1}({\mathbb R})}\|h(x)\|_{L^{2}({\mathbb R})}.
$$
We recall upper bound (\ref{sch}), which yields
$$
\Big\|\int_{0}^{t}e^{(t-s)\{-|p|^{2\alpha}+ibp+a\}}\sqrt{2\pi}p^{2}\widehat{G}(p)
\widehat{f}_{v}(p,s)ds\Big\|_{L^{2}({\mathbb R})}^{2}\leq
$$
$$
e^{2aT}\Big\|\frac{d^{2}G}{dx^{2}}\Big\|_{L^{1}({\mathbb R})}^{2}
\{k\|v(x,s)\|_{L^{2}({\mathbb R}\times [0,T])}\sqrt{T}+
\|h(x)\|_{L^{2}({\mathbb R})}T\}^{2}.
$$
Therefore, we arrive at 
$$
\Big\|\int_{0}^{t}e^{(t-s)\{-|p|^{2\alpha}+ibp+a\}}\sqrt{2\pi}p^{2}\widehat{G}(p)
\widehat{f}_{v}(p,s)ds\Big\|_{L^{2}({\mathbb R}\times [0, T])}^{2}\leq
$$
$$
e^{2aT}\Big\|\frac{d^{2}G}{dx^{2}}\Big\|_{L^{1}({\mathbb R})}^{2}
\{k\|v(x,s)\|_{L^{2}({\mathbb R}\times [0,T])}+
\|h(x)\|_{L^{2}({\mathbb R})}\sqrt{T}\}^{2}T^{2}<\infty
$$
due to the stated assumptions with
$v(x,s)\in W^{1, 2, 2}({\mathbb R}\times [0, T])$. Hence,
\begin{equation}
\label{int0tp2Gpfv}  
\int_{0}^{t}e^{(t-s)\{-|p|^{2\alpha}+ibp+a\}}\sqrt{2\pi}p^{2}\widehat{G}(p)
\widehat{f}_{v}(p,s)ds\in L^{2}({\mathbb R}\times [0, T]).
\end{equation}
By virtue of (\ref{p2u0hpl2}), (\ref{int0tp2Gpfv}) and (\ref{aux2}), we derive
\begin{equation}
\label{up2t12}  
p^{2}\widehat{u}(p,t)\in L^{2}({\mathbb R}\times [0, T]).
\end{equation}  
This proves that
\begin{equation}
\label{d2udx2l2}
\frac{\partial^{2}u}{\partial x^{2}}\in L^{2}({\mathbb R}\times [0, T]).
\end{equation}  
By means of (\ref{aux}),
\begin{equation}
\label{duhdt}
\frac{\partial \widehat{u}}{\partial t}=\{-|p|^{2\alpha}+ibp+a\}\widehat{u}(p,t)+
\sqrt{2\pi}\widehat{G}(p)\widehat{f}_{v}(p,t).
\end{equation}
Clearly, (\ref{upt12}) yields that
\begin{equation}
\label{aupt12}  
a\widehat{u}(p,t)\in L^{2}({\mathbb R}\times [0,T]).
\end{equation}
Let us estimate the norm as
$$
\|ibp\widehat{u}(p,t)\|_{L^{2}({\mathbb R}\times [0,T])}^{2}=b^{2}\int_{0}^{T}
\Big\{\int_{|p|\leq 1}p^{2}|\widehat{u}(p,t)|^{2}dp+
\int_{|p|>1}p^{2}|\widehat{u}(p,t)|^{2}dp\Big\}dt\leq
$$
$$
b^{2}\{\|\widehat{u}(p,t)\|_{L^{2}({\mathbb R}\times [0,T])}^{2}+
\|p^{2}\widehat{u}(p,t)\|_{L^{2}({\mathbb R}\times [0,T])}^{2}\}<\infty
$$
via (\ref{upt12}) and (\ref{up2t12}). Thus,
\begin{equation}
\label{bupt12} 
ibp\widehat{u}(p,t)\in L^{2}({\mathbb R}\times [0,T]).
\end{equation}
Similarly,
$$
\||p|^{2\alpha}\widehat{u}(p,t)\|_{L^{2}({\mathbb R}\times [0,T])}^{2}=
$$
$$
\int_{0}^{T}
\Big\{\int_{|p|\leq 1}|p|^{4\alpha}|\widehat{u}(p,t)|^{2}dp+
\int_{|p|>1}|p|^{4\alpha}|\widehat{u}(p,t)|^{2}dp\Big\}dt\leq
$$
$$
\|\widehat{u}(p,t)\|_{L^{2}({\mathbb R}\times [0,T])}^{2}+
\|p^{2}\widehat{u}(p,t)\|_{L^{2}({\mathbb R}\times [0,T])}^{2}<\infty
$$
due to (\ref{upt12}) and (\ref{up2t12}). Hence,
\begin{equation}
\label{uptalf12} 
|p|^{2\alpha}\widehat{u}(p,t)\in L^{2}({\mathbb R}\times [0,T]).
\end{equation}
Combining (\ref{aupt12}), (\ref{bupt12}) and (\ref{uptalf12}), we arrive at
\begin{equation}
\label{uptalfab12}
(-|p|^{2\alpha}+ibp+a)\widehat{u}(p,t)\in L^{2}({\mathbb R}\times [0,T]).
\end{equation}
Let us turn our attention to obtaing the upper bound on the norm of the
remaining term in the right side of (\ref{duhdt}) using (\ref{fub}) and 
(\ref{lubs}). Obviously,
$$
\|\sqrt{2\pi}\widehat{G}(p)\widehat{f}_{v}(p,t)\|_
{L^{2}({\mathbb R}\times [0,T])}^{2}\leq 2\pi
\|\widehat{G}(p)\|_{L^{\infty}({\mathbb R})}^{2}\int_{0}^{T}
\|F(v(x,t), x)\|_{L^{2}({\mathbb R})}^{2}dt\leq
$$
$$
\|G(x)\|_{L^{1}({\mathbb R})}^{2}\int_{0}^{T}(k\|v(x,t)\|_{L^{2}({\mathbb R})}+
\|h(x)\|_{L^{2}({\mathbb R})})^{2}dt\leq
$$
$$
\|G(x)\|_{L^{1}({\mathbb R})}^{2}\{2k^{2}\|v(x,t)\|_{L^{2}({\mathbb R}\times [0, T])}^{2}+
2\|h(x)\|_{L^{2}({\mathbb R})}^{2}T\}<\infty
$$
under the given conditions with 
$v(x,t)\in W^{1, 2, 2}({\mathbb R}\times [0, T])$. Thus,
\begin{equation}
\label{ghfvhpt}  
\sqrt{2\pi}\widehat{G}(p)\widehat{f}_{v}(p,t)\in L^{2}({\mathbb R}\times [0,T]).
\end{equation}
By virtue of (\ref{duhdt}) along with (\ref{uptalfab12}) and
(\ref{ghfvhpt}),
$$
\frac{\partial \widehat{u}}{\partial t}\in L^{2}({\mathbb R}\times [0,T]),
$$
such that
\begin{equation}
\label{dudtl2}  
\frac{\partial u}{\partial t}\in L^{2}({\mathbb R}\times [0,T])
\end{equation}
as well. Let us recall the definition of the norm (\ref{122n}).
By means of (\ref{uxtl2}), (\ref{d2udx2l2}) and (\ref{dudtl2}), we obtain that
for the function uniquely determined by (\ref{aux}), we have
$$
u(x,t)\in W^{1, 2, 2}({\mathbb R}\times [0, T]).
$$
Therefore, under the stated assumptions problem (\ref{aux}) defines a map
$$
t_{a, b}: W^{1, 2, 2}({\mathbb R}\times [0, T])\to
W^{1, 2, 2}({\mathbb R}\times [0, T]).
$$
Our goal is to demonstrate that under the given conditions this map is a
strict contraction. We choose arbitrarily
$v_{1, 2}(x,t)\in W^{1, 2, 2}({\mathbb R}\times [0, T])$. By virtue of the
argument above,
$u_{1, 2}:=t_{a, b}v_{1, 2}\in W^{1, 2, 2}({\mathbb R}\times [0, T])$.
By means of (\ref{aux}),
$$
\widehat{u_{1}}(p,t)=
$$
\begin{equation}
\label{1aux}
e^{t\{-|p|^{2\alpha}+ibp+a\}}\widehat{u_{0}}(p)+\int_{0}^{t}
e^{(t-s)\{-|p|^{2\alpha}+ibp+a\}}\sqrt{2\pi}\widehat{G}(p)\widehat{f}_{v_{1}}(p,s)ds,
\end{equation}
$$
\widehat{u_{2}}(p,t)=
$$
\begin{equation}
\label{2aux}
e^{t\{-|p|^{2\alpha}+ibp+a\}}\widehat{u_{0}}(p)+\int_{0}^{t}
e^{(t-s)\{-|p|^{2\alpha}+ibp+a\}}\sqrt{2\pi}\widehat{G}(p)\widehat{f}_{v_{2}}(p,s)ds.
\end{equation}
Here $\widehat{f}_{v_{j}}(p,s)$ with $j=1, 2$ stands for the Fourier image of
$F(v_{j}(x,s), x)$ under transform (\ref{ft}). From system
(\ref{1aux}), (\ref{2aux}) we easily derive that
$$
\widehat{u_{1}}(p,t)-\widehat{u_{2}}(p,t)=
$$
\begin{equation}
\label{u1u2hint0t}  
\int_{0}^{t}
e^{(t-s)\{-|p|^{2\alpha}+ibp+a\}}\sqrt{2\pi}\widehat{G}(p)
[\widehat{f}_{v_{1}}(p,s)-\widehat{f}_{v_{2}}(p,s)]ds.
\end{equation}
Let us estimate the norm as
$$
\|\widehat{u_{1}}(p,t)-\widehat{u_{2}}(p,t)\|_{L^{2}({\mathbb R})}\leq
$$
\begin{equation}
\label{int0tetsf12}  
\int_{0}^{t}\|
e^{(t-s)\{-|p|^{2\alpha}+ibp+a\}}\sqrt{2\pi}\widehat{G}(p)
[\widehat{f}_{v_{1}}(p,s)-\widehat{f}_{v_{2}}(p,s)]\|_{L^{2}({\mathbb R})}ds.
\end{equation}
Using inequality (\ref{fub}), we obtain the upper bound as
$$
\|e^{(t-s)\{-|p|^{2\alpha}+ibp+a\}}\sqrt{2\pi}\widehat{G}(p)
[\widehat{f}_{v_{1}}(p,s)-\widehat{f}_{v_{2}}(p,s)]\|_{L^{2}({\mathbb R})}^{2}=
$$
$$
2\pi \int_{-\infty}^{\infty}e^{-2(t-s)|p|^{2\alpha}}e^{2(t-s)a}|\widehat{G}(p)|^{2}
|\widehat{f}_{v_{1}}(p,s)-\widehat{f}_{v_{2}}(p,s)|^{2}dp\leq 
$$
$$
2\pi e^{2aT}\|\widehat{G}(p)\|_{L^{\infty}({\mathbb R})}^{2}\int_{-\infty}^{\infty}
|\widehat{f}_{v_{1}}(p,s)-\widehat{f}_{v_{2}}(p,s)|^{2}dp\leq 
$$
$$
e^{2aT}\|G(x)\|_{L^{1}({\mathbb R})}^{2}\|F(v_{1}(x,s), x)-F(v_{2}(x,s), x)\|_
{L^{2}({\mathbb R})}^{2}.
$$
We recall formula (\ref{lip}). Hence,
\begin{equation}
\label{lipl2}  
\|F(v_{1}(x,s), x)-F(v_{2}(x,s), x)\|_{L^{2}({\mathbb R})}\leq l
\|v_{1}(x,s)-v_{2}(x,s)\|_{L^{2}({\mathbb R})},
\end{equation}
such that
$$
\|e^{(t-s)\{-|p|^{2\alpha}+ibp+a\}}\sqrt{2\pi}\widehat{G}(p)
[\widehat{f}_{v_{1}}(p,s)-\widehat{f}_{v_{2}}(p,s)]\|_{L^{2}({\mathbb R})}\leq
$$
\begin{equation}
\label{eatgvl2}  
e^{aT}l\|G(x)\|_{L^{1}({\mathbb R})}\|v_{1}(x,s)-v_{2}(x,s)\|_{L^{2}({\mathbb R})}.   
\end{equation}
By virtue of (\ref{int0tetsf12}) along with (\ref{eatgvl2}),
$$
\|\widehat{u_{1}}(p,t)-\widehat{u_{2}}(p,t)\|_{L^{2}({\mathbb R})}\leq
$$
$$
e^{aT}l\|G(x)\|_{L^{1}({\mathbb R})}\int_{0}^{T}
\|v_{1}(x,s)-v_{2}(x,s)\|_{L^{2}({\mathbb R})}ds.
$$
The Schwarz inequality yields 
\begin{equation}
\label{sch2}  
\int_{0}^{T}\|v_{1}(x,s)-v_{2}(x,s)\|_{L^{2}({\mathbb R})}ds\leq
\sqrt{\int_{0}^{T}\|v_{1}(x,s)-v_{2}(x,s)\|_{L^{2}({\mathbb R})}^{2}ds}\sqrt{T},
\end{equation}
so that 
$$  
\|\widehat{u_{1}}(p,t)-\widehat{u_{2}}(p,t)\|_{L^{2}({\mathbb R})}\leq
$$
\begin{equation}
\label{u1hu2hl2}  
e^{aT}l\sqrt{T}\|G(x)\|_{L^{1}({\mathbb R})}\|v_{1}(x,t)-v_{2}(x,t)\|_
{L^{2}({\mathbb R}\times [0, T])}.
\end{equation}
Therefore,
$$
\|u_{1}(x,t)-u_{2}(x,t)\|_{L^{2}({\mathbb R}\times [0, T])}^{2}=\int_{0}^{T}
\|\widehat{u_{1}}(p,t)-\widehat{u_{2}}(p,t)\|_{L^{2}({\mathbb R})}^{2}dt\leq
$$
\begin{equation}
\label{u1u2l2v1v2}
e^{2aT}l^{2}T^{2}\|G(x)\|_{L^{1}({\mathbb R})}^{2}
\|v_{1}(x,t)-v_{2}(x,t)\|_{L^{2}({\mathbb R}\times [0, T])}^{2}.
\end{equation}  
By means of (\ref{u1u2hint0t}),
$$
p^{2}[\widehat{u_{1}}(p,t)-\widehat{u_{2}}(p,t)]=
\int_{0}^{t}
e^{(t-s)\{-|p|^{2\alpha}+ibp+a\}}\sqrt{2\pi}p^{2}\widehat{G}(p)
[\widehat{f}_{v_{1}}(p,s)-\widehat{f}_{v_{2}}(p,s)]ds.
$$
We obtain the upper bound on the norm as
$$
\|p^{2}[\widehat{u_{1}}(p,t)-\widehat{u_{2}}(p,t)]\|_{L^{2}({\mathbb R})}\leq
$$
\begin{equation}
\label{int0tetsf122}  
\int_{0}^{t}\|
e^{(t-s)\{-|p|^{2\alpha}+ibp+a\}}\sqrt{2\pi}p^{2}\widehat{G}(p)
[\widehat{f}_{v_{1}}(p,s)-\widehat{f}_{v_{2}}(p,s)]\|_{L^{2}({\mathbb R})}ds.
\end{equation}
Let us use inequality (\ref{fub1}) to derive the estimate from above
$$
\|e^{(t-s)\{-|p|^{2\alpha}+ibp+a\}}\sqrt{2\pi}p^{2}\widehat{G}(p)
[\widehat{f}_{v_{1}}(p,s)-\widehat{f}_{v_{2}}(p,s)]\|_{L^{2}({\mathbb R})}^{2}=
$$
$$
2\pi \int_{-\infty}^{\infty}e^{-2(t-s)|p|^{2\alpha}}e^{2(t-s)a}|p^{2}\widehat{G}(p)|^{2}
|\widehat{f}_{v_{1}}(p,s)-\widehat{f}_{v_{2}}(p,s)|^{2}dp\leq 
$$
$$
2\pi e^{2aT}\|p^{2}\widehat{G}(p)\|_{L^{\infty}({\mathbb R})}^{2}\int_{-\infty}^{\infty}
|\widehat{f}_{v_{1}}(p,s)-\widehat{f}_{v_{2}}(p,s)|^{2}dp\leq 
$$
$$
e^{2aT}\Big\|\frac{d^{2}G}{dx^{2}}\Big\|_{L^{1}({\mathbb R})}^{2}
\|F(v_{1}(x,s), x)-F(v_{2}(x,s), x)\|_{L^{2}({\mathbb R})}^{2}.
$$
By virtue of formula (\ref{lipl2}),
$$
\|e^{(t-s)\{-|p|^{2\alpha}+ibp+a\}}\sqrt{2\pi}p^{2}\widehat{G}(p)
[\widehat{f}_{v_{1}}(p,s)-\widehat{f}_{v_{2}}(p,s)]\|_{L^{2}({\mathbb R})}\leq 
$$
\begin{equation}
\label{etsp2gpeat}
e^{aT}l\Big\|\frac{d^{2}G}{dx^{2}}\Big\|_{L^{1}({\mathbb R})}\|v_{1}(x,s)-v_{2}(x,s)\|_
{L^{2}({\mathbb R})}.
\end{equation}
Using (\ref{int0tetsf122}) along with (\ref{sch2}) and (\ref{etsp2gpeat}), we
derive that
$$
\|p^{2}[\widehat{u_{1}}(p,t)-\widehat{u_{2}}(p,t)]\|_{L^{2}({\mathbb R})}\leq
$$
\begin{equation}
\label{p2u1hu2hpt}  
e^{aT}\sqrt{T}l\Big\|\frac{d^{2}G}{dx^{2}}\Big\|_{L^{1}({\mathbb R})}
\|v_{1}(x,t)-v_{2}(x,t)\|_{L^{2}({\mathbb R}\times [0, T])},
\end{equation}
such that
$$
\Big\|\frac{\partial^{2}}{\partial x^{2}}[u_{1}(x,t)-u_{2}(x,t)]\Big\|_
{L^{2}({\mathbb R}\times [0, T])}^{2}=\int_{0}^{T}
\|p^{2}[\widehat{u_{1}}(p,t)-\widehat{u_{2}}(p,t)]\|_{L^{2}({\mathbb R})}^{2}dt\leq     $$
\begin{equation}
\label{u1u2l2v1v22}
e^{2aT}l^{2}T^{2}\Big\|\frac{d^{2}G}{dx^{2}}\Big\|_{L^{1}({\mathbb R})}^{2}
\|v_{1}(x,t)-v_{2}(x,t)\|_{L^{2}({\mathbb R}\times [0, T])}^{2}.
\end{equation}  
From (\ref{u1u2hint0t}) we easily deduce that
$$
\frac{\partial}{\partial t}[\widehat{u_{1}}(p,t)-\widehat{u_{2}}(p,t)]=
$$
$$
\{-|p|^{2\alpha}+ibp+a\}[\widehat{u_{1}}(p,t)-\widehat{u_{2}}(p,t)]+
\sqrt{2\pi}\widehat{G}(p)[\widehat{f}_{v_{1}}(p,t)-\widehat{f}_{v_{2}}(p,t)],
$$
so that
$$
\Big\|\frac{\partial}{\partial t}[\widehat{u_{1}}(p,t)-\widehat{u_{2}}(p,t)]
\Big\|_{L^{2}({\mathbb R})}\leq a\|\widehat{u_{1}}(p,t)-\widehat{u_{2}}(p,t)\|_
{L^{2}({\mathbb R})}+
$$
$$      
|b|\|p[\widehat{u_{1}}(p,t)-\widehat{u_{2}}(p,t)]\|_
{L^{2}({\mathbb R})}+
\||p|^{2\alpha}[\widehat{u_{1}}(p,t)-\widehat{u_{2}}(p,t)]\|_{L^{2}({\mathbb R})}+    
$$
\begin{equation}
\label{ddtu1hptu2hpt}
\sqrt{2\pi}\|\widehat{G}(p)[\widehat{f}_{v_{1}}(p,t)-\widehat{f}_{v_{2}}(p,t)]\|_
{L^{2}({\mathbb R})}.
\end{equation}  
By means of (\ref{u1hu2hl2}), the first term in the right side of
(\ref{ddtu1hptu2hpt}) can be bounded from above by
\begin{equation}
\label{au1hu2hl2}  
aQe^{aT}\sqrt{T}l\|v_{1}(x,t)-v_{2}(x,t)\|_{L^{2}({\mathbb R}\times [0, T])},
\end{equation}
where $Q$ is defined in (\ref{q}).

\noindent
Let us estimate the norm as
$$
\|p[\widehat{u_{1}}(p,t)-\widehat{u_{2}}(p,t)]\|_{L^{2}({\mathbb R})}^{2}=
$$
$$
\int_{|p|\leq 1}p^{2}|\widehat{u_{1}}(p,t)-\widehat{u_{2}}(p,t)|^{2}dp+
\int_{|p|>1}p^{2}|\widehat{u_{1}}(p,t)-\widehat{u_{2}}(p,t)|^{2}dp\leq 
$$
$$  
\|\widehat{u_{1}}(p,t)-\widehat{u_{2}}(p,t)\|_{L^{2}({\mathbb R})}^{2}+
\|p^{2}[\widehat{u_{1}}(p,t)-\widehat{u_{2}}(p,t)]\|_{L^{2}({\mathbb R})}^{2}.
$$
We recall inequalities (\ref{u1hu2hl2}) and  (\ref{p2u1hu2hpt}). They
yield the upper bound on the second term in the right side of
(\ref{ddtu1hptu2hpt}) given by
\begin{equation}
\label{bu1hu2hl2}  
|b|Qe^{aT}\sqrt{T}l\|v_{1}(x,t)-v_{2}(x,t)\|_{L^{2}({\mathbb R}\times [0, T])}.
\end{equation}
Similarly,
$$
\||p|^{2\alpha}[\widehat{u_{1}}(p,t)-\widehat{u_{2}}(p,t)]\|_{L^{2}({\mathbb R})}^{2}=
$$
$$
\int_{|p|\leq 1}|p|^{4\alpha}|\widehat{u_{1}}(p,t)-\widehat{u_{2}}(p,t)|^{2}dp+
\int_{|p|>1}|p|^{4\alpha}|\widehat{u_{1}}(p,t)-\widehat{u_{2}}(p,t)|^{2}dp\leq 
$$
$$  
\|\widehat{u_{1}}(p,t)-\widehat{u_{2}}(p,t)\|_{L^{2}({\mathbb R})}^{2}+
\|p^{2}[\widehat{u_{1}}(p,t)-\widehat{u_{2}}(p,t)]\|_{L^{2}({\mathbb R})}^{2}.
$$
By virtue of (\ref{u1hu2hl2}) and  (\ref{p2u1hu2hpt}), the third term in the
right side of (\ref{ddtu1hptu2hpt}) can be estimated from above by
\begin{equation}
\label{bu1hu2hl3}  
Qe^{aT}\sqrt{T}l\|v_{1}(x,t)-v_{2}(x,t)\|_{L^{2}({\mathbb R}\times [0, T])}.
\end{equation}
Using (\ref{fub}) along with (\ref{lipl2}), we obtain that
$$
2\pi\int_{-\infty}^{\infty}|\widehat{G}(p)|^{2}
|\widehat{f}_{v_{1}}(p,t)-\widehat{f}_{v_{2}}(p,t)|^{2}dp\leq
$$
$$
2\pi
\|\widehat{G}(p)\|_{L^{\infty}({\mathbb R})}^{2}\int_{-\infty}^{\infty}
|\widehat{f}_{v_{1}}(p,t)-\widehat{f}_{v_{2}}(p,t)|^{2}dp\leq
$$
$$
\|G(x)\|_{L^{1}({\mathbb R})}^{2}\|F(v_{1}(x,t), x)-F(v_{2}(x,t), x)\|_
{L^{2}({\mathbb R})}^{2}\leq
$$
$$
\|G(x)\|_{L^{1}({\mathbb R})}^{2}l^{2}\|v_{1}(x,t)-v_{2}(x,t)\|_{L^{2}({\mathbb R})}^{2}.  
$$
Hence, the fourth term in the right side of (\ref{ddtu1hptu2hpt}) can be bounded
from above by
\begin{equation}
\label{bu1hu2hl4}  
Ql\|v_{1}(x,t)-v_{2}(x,t)\|_{L^{2}({\mathbb R})}.
\end{equation}
Combining (\ref{au1hu2hl2}), (\ref{bu1hu2hl2}), (\ref{bu1hu2hl3}) and
(\ref{bu1hu2hl4}), we arrive at
$$
\Big\|\frac{\partial}{\partial t}[\widehat{u_{1}}(p,t)-\widehat{u_{2}}(p,t)]
\Big\|_{L^{2}({\mathbb R})}\leq
$$
$$
Qe^{aT}\sqrt{T}l\{a+|b|+1\}\|v_{1}(x,t)-v_{2}(x,t)\|_
{L^{2}({\mathbb R}\times [0, T])}+
Ql\|v_{1}(x,t)-v_{2}(x,t)\|_{L^{2}({\mathbb R})}.
$$
This enables us to estimate the norm as
$$
\Big\|\frac{\partial}{\partial t}(u_{1}(x,t)-u_{2}(x,t))\Big\|_
{L^{2}({\mathbb R}\times [0, T])}^{2}=\int_{0}^{T}
\Big\|\frac{\partial}{\partial t}[\widehat{u_{1}}(p,t)-\widehat{u_{2}}(p,t)]
\Big\|_{L^{2}({\mathbb R})}^{2}dt\leq           
$$
\begin{equation}
\label{ddtu1u2l2}  
2Q^{2}l^{2}[e^{2aT}T^{2}\{a+|b|+1\}^{2}+1]
\|v_{1}(x,t)-v_{2}(x,t)\|_{L^{2}({\mathbb R}\times [0, T])}^{2}.
\end{equation}
Let us recall the definition of the norm (\ref{122n}) and combine upper
bounds (\ref{u1u2l2v1v2}), (\ref{u1u2l2v1v22}) and (\ref{ddtu1u2l2}). 
A straightforward computation yields that
$$
\|u_{1}-u_{2}\|_{W^{1, 2, 2}({\mathbb R}\times [0, T])}\leq
$$
\begin{equation}
\label{contr}
Ql\sqrt{T^{2}e^{2aT}(1+2[a+|b|+1]^{2})+1}
\|v_{1}-v_{2}\|_{W^{1, 2, 2}({\mathbb R}\times [0, T])}.
\end{equation}
The constant in the right side of inequality (\ref{contr}) is less than one
according to (\ref{qlt}). Therefore, under the given conditions equation
(\ref{aux}) defines the map
$$
t_{a, b}: W^{1, 2, 2}({\mathbb R}\times [0, T])\to W^{1, 2, 2}
({\mathbb R}\times [0, T]),
$$
which is a strict contraction. Its unique fixed
point $w(x,t)$ is the only solution of
problem (\ref{h}), (\ref{ic}) in
$W^{1, 2, 2}({\mathbb R}\times [0, T])$. \hfill\lanbox

\bigskip

\noindent
{\it Proof of Corollary 1.4.} The validity of the statement of the Corollary
follows from the fact that the constant in the right side of estimate
(\ref{contr}) is independent of the initial condition (\ref{ic}), such that
problem (\ref{h}), (\ref{ic}) has a unique solution
$w(x,t)\in W^{1, 2, 2}({\mathbb R}\times {\mathbb R}^{+})$. Let us suppose
that $w(x,t)\equiv 0$ for $x\in {\mathbb R}$ and $t\in {\mathbb R}^{+}$.
This will imply the contradiction to our assumption that
$\hbox{supp}\widehat{F(0, x)}\cap \hbox{supp}\widehat{G}$ is a set of
nonzero Lebesgue measure on the real line. \hfill\lanbox

\bigskip


\centerline{\bf 3. Acknowledgement}

\bigskip

\noindent
V.V. is grateful to Israel Michael Sigal for the partial support by the
NSERC grant NA 7901.

\bigskip


\centerline{\bf Conflicts of interest}

\bigskip

\noindent
This work does not have any conflicts of interest.


\bigskip

\begin{thebibliography}{99}

\bibitem{ABVV10}
N. ~Apreutesei, N. ~Bessonov, V. ~Volpert, V. ~Vougalter.
{\em Spatial structures and generalized travelling waves for an integro-   
differential equation}, Discrete Contin. Dyn. Syst. Ser. B,
{\bf 13} (2010), no. 3, 537--557.

\bibitem{BNPR09}
H. ~Berestycki, G. ~Nadin, B. ~Perthame, L. ~Ryzhik.
{\em The non-local Fisher-KPP equation: travelling waves and steady states},
Nonlinearity, {\bf 22} (2009), no. 12, 2813--2844.

\bibitem{BHN05}
H. ~Berestycki, F. ~Hamel, N. ~Nadirashvili.
{\em The speed of propagation for KPP type problems. I: Periodic 
framework}, J. Eur. Math. Soc. (JEMS), {\bf 7} (2005), no. 2, 173--213. 


\bibitem{BO86}
H. ~Brezis, L. ~Oswald. {\em Remarks on sublinear elliptic equations},
Nonlinear Anal., {\bf 10} (1986), no. 1, 55--64.  

\bibitem{2}
 B. ~Carreras, V. ~Lynch, G. ~Zaslavsky. {\em Anomalous diffusion and exit time
distribution of particle tracers in plasma turbulence model}, Phys. Plasmas,
{\bf 8} (2001), 5096--5103.


\bibitem{DMV05} A. ~Ducrot, M. ~Marion, V. ~Volpert.  {\em Syst\'emes de
r\'eaction-diffusion sans propri\'et\'e de Fredholm},
C. R. Math. Acad. Sci. Paris,  {\bf 340} (2005), no. 9, 659--664.


\bibitem{DMV08} A. ~Ducrot, M. ~Marion, V. ~Volpert.
{\em Reaction-diffusion problems with non Fredholm operators},
Adv. Differ. Equations, {\bf 13} (2008), no. 11-12, 1151--1192.


\bibitem{EV20} M. ~Efendiev, V. ~Vougalter.
{\em Solvability of some integro-differential equations with drift},
Osaka J. Math., {\bf 57} (2020), no. 2, 247--265.


\bibitem{EV22} M. ~Efendiev, V. ~Vougalter.
{\em Solvability of Some Integro-Differential Equations with Drift and
Superdiffusion}, J. Dynam. Differential Equations, {\bf 36} (2024), no. 1,
353--373.


\bibitem{K64} M.A. ~Krasnosel'skii.
{\em Topological methods in the theory of nonlinear integral equations}.
International Series of Monographs on Pure and Applied Mathematics.
Pergamon Press, {\bf XI}, (1964), 395 pp.


\bibitem{LL97} E.H. ~Lieb, M.~Loss.
{\em Analysis. Grad. Stud. Math.}, {\bf 14}, American
 Mathematical Society, Providence, RI (1997), 278 pp.


\bibitem{3}
P. Manandhar, J. Jang, G.C. Schatz, M.A. Ratner, S. Hong. {\em Anomalous surface
diffusion in nanoscale direct deposition processes}, Phys. Rev. Lett., {\bf 90}
(2003), 4043--4052.


\bibitem{MK00}
R. ~Metzler, J. ~Klafter. {\em The random walk's guide to anomalous diffusion:
a fractional dynamics approach}, Phys. Rep., {\bf 339} (2000), 1--77.

\bibitem{4}
J. Sancho, A. Lacasta, K. Lindenberg, I. Sokolov, A. Romero. {\em Diffusion on
a solid surface: anomalous is normal}, Phys. Rev. Lett., {\bf 92} (2004),
250601.

\bibitem{5}
H. ~Scher, E. ~Montroll. {\em Anomalous transit-time dispersion in amorphous
solids}, Phys. Rev. B, {\bf 12} (1975), 2455--2477.

\bibitem{1}
 T. Solomon, E. Weeks, H. Swinney. {\em Observation of anomalous diffusion and
L\'evy flights in a two-dimensional rotating flow}, Phys. Rev. Lett., {\bf 71}
(1993), 3975--3978.

\bibitem{VV130} V. ~Volpert, V. ~Vougalter. {\em Emergence and propagation
of patterns in nonlocal reaction-diffusion equations arising in the theory
of speciation.} Dispersal, individual movement and spatial ecology,
Lecture Notes in Math., {\bf 2071} (2013), Springer, Heidelberg, 331--353.

\bibitem{VNN10} V.A. ~Volpert, Y. ~Nec, A.A. ~Nepomnyashchy.
{\em Exact solutions in front propagation problems with superdiffusion}, Phys.
D, {\bf 239} (2010), no. 3--4, 134--144.

\bibitem{VNN13} V.A. ~Volpert, Y. ~Nec, A.A. ~Nepomnyashchy. {\em Fronts in
anomalous diffusion-reaction systems}, Philos. Trans. R. Soc. Lond. Ser. A
Math. Phys. Eng. Sci., {\bf 371} (2013), no. 1982, 20120179, 18 pp.


\bibitem{VV18} V. ~Vougalter, V. ~Volpert. {\em Solvability of some integro-
differential equations with anomalous diffusion.} Regularity and stochasticity
of nonlinear dynamical systems, Nonlinear Syst. Complex., {\bf 21},
Springer, Cham (2018), 1--17.

\bibitem{VV21} V. ~Vougalter, V. ~Volpert. {\em Solvability of some integro-
differential equations with anomalous diffusion and transport},
Anal. Math. Phys., {\bf 11} (2021), no. 3, Paper No. 135, 26 pp.     

  
\end{thebibliography}
\end{document}